\documentclass{emsprocart}
\usepackage{hyperref}
\usepackage{cite}
\usepackage{graphics}
\usepackage{color}

\newcommand{\e}{^\varepsilon}
\newcommand{\eps}{{\varepsilon}}
\newcommand{\ds}{\displaystyle}
\newcommand{\I}{\mathcal{I}\e}

\renewcommand{\a}{\alpha}
\renewcommand{\b}{\beta}

\newcommand{\cupl}{\bigcup\limits}

\newcommand{\suml}{\sum\limits}

\newcommand{\liml}{\lim\limits}

\newcommand{\Om}{\widetilde{\Omega}}
\newcommand{\Ga}{\widetilde{\Gamma}}
 
\renewcommand{\phi}{\varphi}

\newcommand{\x}{\textbf{x}}

\renewcommand{\d}{\hspace{1pt}\mathrm{d}}
 
\contact[giuseppe.cardone@unisannio.it]{Giuseppe Cardone, Department of Engineering, University of Sannio, Corso Garibaldi 107, 82100 Benevento, Italy}

\contact[andrii.khrabustovskyi@kit.edu]{Andrii Khrabustovskyi, Institute of Analysis, Karlsruhe Institute of Technology, Englerstra{\ss}e. 2, 76131 Karlsruhe , Germany}

\newtheorem{theorem}{Theorem}[section]
\newtheorem{corollary}[theorem]{Corollary}
\newtheorem{lemma}[theorem]{Lemma}

\theoremstyle{definition}

\newtheorem{remark}[theorem]{Remark}

\title[Example of periodic Neumann waveguide with gap in spectrum]{Example of periodic Neumann waveguide with gap in spectrum}

\author[G.~Cardone, A.~Khrabustovskyi]{Giuseppe Cardone, Andrii Khrabustovskyi\thanks{{G.C. is a member of GNAMPA (INdAM).} A.K. is supported by the Deutsche
Forschungsgemeinschaft (DFG) through CRC 1173.}}

\begin{document}

\begin{abstract}
In this note we investigate spectral properties of a periodic waveguide $\Omega\e$ ($\eps$ is a small parameter) obtained from a straight strip 
by attaching an array of $\eps$-periodically distributed identical protuberances having ``room-and-passage'' geometry. In the current work we consider the operator $\mathcal{A}\e=-\rho\e\Delta_{\Omega\e}$, where $\Delta_{\Omega\e}$ is the Neumann Laplacian in $\Omega\e$, the weight $\rho\e$ is equal to $1$ everywhere except  the union of the ``rooms''. We will prove that the spectrum of $\mathcal{A}\e$ has at least one gap as $\eps$ is small enough provided certain conditions on the weight $\rho\e$ and the sizes of attached protuberances hold.
\end{abstract}



\maketitle

\begin{flushright}\it\small
This paper is dedicated to Pavel Exner on the occasion of his jubilee 
\end{flushright}

\section{Introduction}

It is a well-known fact (see, e.g., \cite{Kuchment}) that the spectrum of 
periodic elliptic self-adjoint differential operators has band structure, 
i.e. it is a locally finite union of compact intervals called \textit{bands}. In general the bands may overlap, otherwise we have a \textit{gap} in the spectrum - a bounded open interval having an empty intersection with the spectrum, but with ends belonging to it.

The presence of gaps in the spectrum is not guaranteed. For instance, the spectrum of the Laplace operator in $L^2(\mathbb{R}^n)$ has no gaps: $\sigma(-\Delta_{\mathbb{R}^n})=[0,\infty)$. Therefore an interesting question arises here: to construct examples of periodic operators with non-void spectral gaps. This question is motivated by various applications, since the presence of gaps is important for the description of wave processes which are governed by differential operators
under consideration: if the wave frequency belongs to a gap then the corresponding wave cannot propagate in the medium. This feature is a main requirement for so-called photonic crystals, which are materials with periodic dielectric structure {extensively} investigating in recent years.

The problem of existence of spectral gaps for various periodic operators has been actively studied since mid 90th. We refer to the overview \cite{HempelPost}, where one can find a lot of examples and references around this topic. 

In the last years there appeared many works, where the problem of opening of spectral gaps for  operators posed in unbounded domains with a waveguide geometry (strips, tubes, {graph-like} domains, etc.) is studied, see, e.g., \cite{BNR,Borisov2,Cardone1,Cardone2,EP,Nazarov1+,Nazarov2,Pankrashkin,Yoshi}. 
The studies of physical processes (e.g., quantum particle motion) in such domains are of a great physical and mathematical interest because of the big progress in microelectronics during the last decade. We refer to the recent monograph \cite{EK} concerning spectral properties of quantum waveguides.

The simplest way to open up a gap is either to perturb a straight cylinder by a periodic nucleation of small voids (or making other ``small'' perturbation) \cite{BNR,Nazarov2} or to consider a waveguide consisting of an array of identical compact domains connected by narrow ``bridges'' \cite{Nazarov1+,Pankrashkin}. In the first case one has small gaps separating large bands, in the second case one gets large gaps and small bands.

In the current paper we present another example of Neumann waveguide with a gap in the spectrum; the geometry of this waveguide essentially differs from previously studied examples. We are motivated by our recent work \cite{CarKhrab}, where the spectrum of some Neumann problem was studied in a \textit{bounded} domain perturbed by a lot of identical protuberances each of them consisting of two subsets - ``room'' and ``passage'' (in the simplest case, ``room'' is a small square and ``passage'' is a narrow rectangle connecting the ``room'' with the main domain). Peculiar spectral properties of so perturbed domains were observed for the first time by R.~Courant and D.~Hilbert \cite{CH}. Domains with "room-and-passage"-like geometry are widely used in order to construct examples illustrating various phenomena in Sobolev spaces theory and in spectral theory (see, for example, \cite{Fraenkel,HSS}).

\begin{figure}[t]
\begin{picture}(300,85)      
     \scalebox{0.45}{\includegraphics{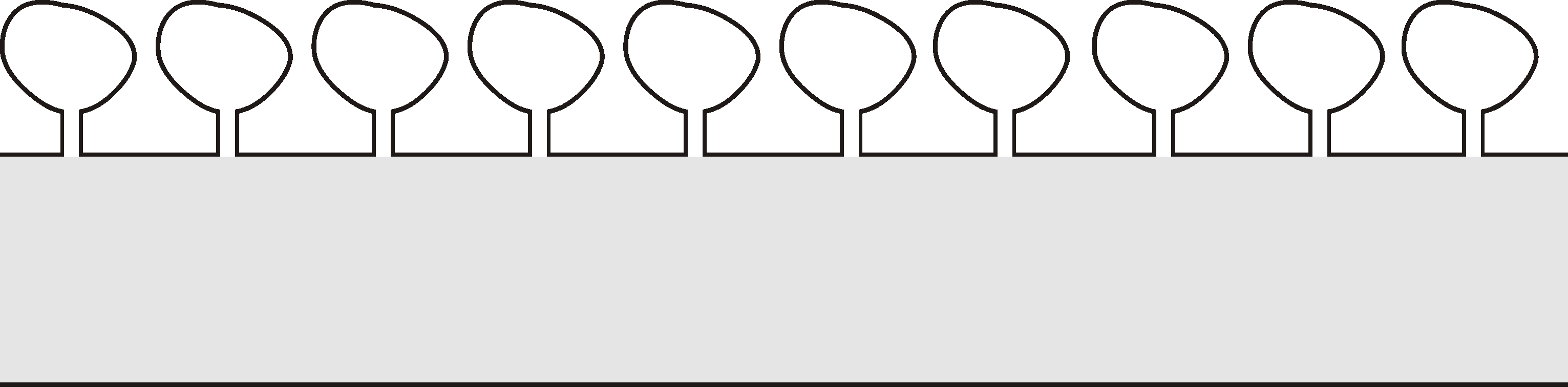}}
     \put(-40, 20){$\Omega$}

     \put(-145, 53){$B_i\e$}
     \put(-133,56){\vector(1,1){10}}

     \put(-104, 53){$T_i\e$}
     \put(-105,55){\vector(-1,0){15}}
     
     \put(-330,30){\vector(0,1){20}}
     \put(-330,20){\vector(0,-1){20}}     
     \put(-333, 21){$L$}
     
     \put(-297,45){\vector(1,0){14}}
     \put(-304,45){\vector(-1,0){14}}     
     \put(-303, 42){$\eps$}

\end{picture}
\caption{\label{figure}The waveguide $\Omega\e$}
\end{figure}
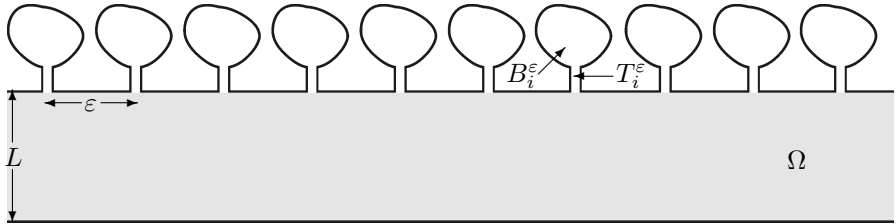

Our goal is to show that perturbing a straight strip by a periodic array of ``room-and-passage'' protuberances one may open a spectral gap. Namely, we consider a strip of a width $L>0$ and perturb it by 
a family of small identical protuberances, $\eps$-periodically distributed along the strip axis. Here $\eps>0$ is a small parameter. Each protuberance has ``room-and-passage'' geometry. We denote the obtained domain by $\Omega\e$ {(see Figure \ref{figure})}. In $\Omega\e$ we consider the operator $\mathcal{A}\e=-{\rho\e}\Delta_{\Omega\e}$, where $\Delta_{\Omega\e}$ is the Neumann Laplacian in $L^2(\Omega\e)$. The weight $\rho\e$ is equal to $1$ everywhere except  the union of the ``rooms'', where it is equal to the constant $\varrho\e>0$.

The main result: we will prove that under suitable assumptions on $L,\ \varrho\e$ and sizes of ``rooms'' and ``passages'' the spectrum of $\mathcal{A}\e$ converges 
to the spectrum of a certain spectral problem on the initial strip containing the spectral parameter in boundary conditions. Its spectrum has the form $[0,\infty)\setminus (\alpha,\beta)$, where
$(\alpha,\beta)$ is a non-empty bounded interval. This, in particular, implies  
at least one gap in the spectrum of $\mathcal{A}\e$ for small enough $\eps$.

\section{Setting of the problem and main result}

In what follows by $x$ and $\mathbf{x}=(x_1,x_2)$ we denote the Cartesian coordinates in $\mathbb{R}$ and $\mathbb{R}^2$, correspondingly. 

By $\eps$ we denote a small parameter. To simplify the proof of the main theorem we suppose that it takes values from the discrete set $\mathcal{E}=\left\{\eps:\ \eps^{-1}\in\mathbb{N}\right\}$. 
The general case needs slight modifications. 

We consider the unbounded strip $\Omega\subset\mathbb{R}^2$ of the width $L>0$:
\begin{align*}
\Omega=\left\{\mathbf{x}\in\mathbb{R}^2: -L<x_2<0\right\}.
\end{align*}
By $\Gamma$ we denote its upper boundary: $\Gamma=\left\{\mathbf{x}\in\mathbb{R}^2: x_2=0\right\}.$

Let $b\e$, $d\e$, $h\e$ be positive constants,
$B$ be an open bounded domain in $\mathbb{R}^{2}$ having Lipschitz boundary and satisfying
\begin{gather}\label{ass11}
B\subset \left\{\x\in\mathbb{R}^2:\ x_1\in\left(-{1/2},{1/2}\right),\ x_2>0\right\},\\\label{ass12}
\exists R\in \left(0,1\right):\ \left\{\x\in\mathbb{R}^2:\ x_1\in\left(-{R/2},{R/2}\right),\ x_2=0\right\}\subset \partial B,\\
\label{ass14}
R^{-1}d\e\leq b\e\leq \eps,\\
\label{ass15}
h\e\to 0\text{ as }\eps\to 0.
\end{gather}

For $i\in\mathbb{Z}$ we set:
\begin{itemize}

\item[] $B_i\e=\left\{\x\in\mathbb{R}^2:\ \ds{1\over b\e}\left(\x- \widetilde{\x}^{i,\eps}\right)\in B,\text{ where } \widetilde{\x}^{i,\eps}=(i\eps+\eps/2,h\e)\right\}$\quad (``room''),

\item[] $T_i\e=\left\{\x\in\mathbb{R}^2:\ \ds |x_1-i\eps-\eps/2|<{d\e\over 2},\ 0\leq x_2\leq h\e\right\}$ \quad (``passage'').

\end{itemize} 
Conditions \eqref{ass11}-\eqref{ass14} imply that the "rooms" are pairwise disjoint and guarantee correct gluing of the $i$-th "room" and the $i$-th "passage" (the upper face of $T_i\e$ is contained in $\partial B_i\e$).  Moreover, the distance between the neighbouring "passages" is not too small, namely
for $i\not=j$ one has $\mathrm{dist}( {T_i\e,T_j\e})\geq \eps- 
d\e\geq \eps\left (1-R\right)$.

Attaching the "rooms" and "passages" to $\Omega$ we obtain the perturbed 
domain 
$$\Omega\e=\Omega\cup \left({\cupl_{i\in \mathbb{Z}}\left(T_i\e\cup B_i\e\right)}\right).$$

Let us define accurately the operator $\mathcal{A}\e$. We denote by ${H}\e$ the Hilbert space of functions from $L^2(\Omega\e)$ endowed with a scalar product
\begin{gather}\label{He}
(u,v)_{{H}\e}=\int_{\Omega\e} u(\x)\overline{v(\x)}  (\rho \e(\x))^{-1} \d \x,
\end{gather}
where the function $\rho\e(\x)$ is defined as follows:
$$\rho\e(\x)=\begin{cases}1,&\x\in\Omega\cup\left(\cupl_{i\in\mathbb{Z}}T_i\e\right),\\\varrho\e,& \x\in \cupl_{i\in\mathbb{Z}}B_i\e,\end{cases}\quad \varrho\e>0\text{ is a constant.}
$$
By $\mathfrak{a}\e$ we denote the sesquilinear form in ${H}\e$ defined by  
\begin{gather}\label{ae}
\mathfrak{a}\e[u,v]=\int_{\Omega\e} \nabla u\cdot\overline{\nabla  v} \d \x,\quad \mathrm{dom}(\mathfrak{a}\e)=H^1(\Omega\e).
\end{gather}
The form  $\mathfrak{a}\e$ is densely defined, closed, positive and symmetric. We denote by
$\mathcal{A}\e$ the operator associated with this form, i.e.
\begin{gather*}
(\mathcal{A}\e u,v)_{{H}\e}=
\mathfrak{a}\e[u,v],\quad\forall u\in
\mathrm{dom}(\mathcal{A}\e),\ \forall v\in
\mathrm{dom}(\mathfrak{a}\e).
\end{gather*}
In other words, the operator $\mathcal{A}\e$ is defined by the operation $-{\rho\e}\Delta$ in $\Omega\e$ and
the  Neumann boundary conditions on $\partial\Omega\e$.

The goal of this work is to describe the behaviour of the spectrum $\sigma(\mathcal{A}\e)$ as $\eps\to 0$  
under the assumption that the following limits exist and are positive:
\begin{gather}\label{qere+}
\a:=\liml_{\eps\to 0}{d\e \varrho\e\over h\e (b\e)^2 |B| },\quad r:=\liml_{\eps\to 0}{(b\e)^2 |B|\over  \eps\varrho\e },\quad \a>0,\ r>0.
\end{gather}
Also it is supposed that $d\e$ tends to zero not very fast, namely
$\liml_{\eps\to 0}\eps \ln d\e= 0.$ The meaning of this condition and the meaning of $\a$ and $r$ are explained in \cite{CarKhrab}.

Now, we introduce the limit operator. By ${H}$ we denote the Hilbert space of functions from $L^2(\Omega)\oplus L^2(\Gamma)$ endowed with the scalar product
\begin{gather}\label{Hlim}
(U,V)_{{H}}=\int_{\Omega}u_1(\x)\overline{v_1(\x)} \d \x+\int_\Gamma u_2(x)\overline{v_2(x)} r \d x,\ U=\left(u_1,u_2\right),V=(v_1,v_2).
\end{gather}
We introduce the sesquilinear form $\mathfrak{a}$ in ${H}$ by 
\begin{align}\label{alim}
\mathfrak{a}[U,V]=\int_\Omega \nabla u_1\cdot \overline{\nabla v_1}\d 
\x+\int_\Gamma \a r\left(u_1|_\Gamma-u_2\right)\overline{\left(v_1|_\Gamma-v_2\right)}\d x
\end{align}
with $\mathrm{dom}(\mathfrak{a} )=H^1(\Omega)\oplus L^2(\Gamma)$.
Here by $u|_\Gamma$ we denote the trace of  $u$ on $\Gamma$.
We denote by $\mathcal{A}$ the self-adjoint operator associated with this form. 

Formally, the eigenvalue equation $\mathcal{A} U=\lambda U$ can be written as follows:
\begin{gather*}
\begin{cases}
-\Delta u_1=\lambda u_1&\text{ in }\Omega,\\ 
\ds{\partial u_1\over\partial  n}=\a r(u_2-u_1)&\text{ on }\Gamma,\\
\a (u_2-u_1)=\lambda u_2&\text{ on }\Gamma,\\
\ds{\partial u_1\over\partial n}=0&\text{ on }\partial\Omega\setminus\Gamma.
\end{cases}
\end{gather*}
Here $n$ is the outward-pointing unit normal.

\begin{remark}
Spectral properties of so defined operators $\mathcal{A}$ were investigated in \cite{CarKhrab,KhrabPlum}. In \cite{CarKhrab} one considered the case of a bounded domain $\Omega$, $\Gamma$ is a flat subset of $\partial\Omega$. In this case the discrete spectrum of $\mathcal{A}$ consists of two sequences; one sequence accumulates at $\infty$, while the other one converges to $\a$, which is the only point of the essential spectrum.

In \cite{KhrabPlum} one considered\footnote{In fact, in \cite{KhrabPlum} the Dirichlet conditions on $\partial\Omega$ are prescribed, but similar results can be easily obtained for the Neumann conditions too -- cf. \cite[Remark 3.2]{KhrabPlum}.}, in particular, the case, when $\Omega$ is a straight unbounded strip, the line $\Gamma$ is parallel to its axis and divides $\Omega$ on two unbounded strips. In this case the spectrum of $\mathcal{A}$ turns out to be a union of the interval $[0,\a]$ and the ray $[\b,\infty)$, where  $\b>\a$ provided $\a<\left(\pi\over L-L_\Gamma\right)^2$. Here $L$ is strip width and $L_\Gamma\in (0,L)$ is a distance from $\Gamma$ to $\partial\Omega$. 
\end{remark}

Using the same arguments as in \cite{KhrabPlum} we arrive at the following formula for the spectrum of \textit{our} operator $\mathcal{A}$:
\begin{gather}\label{aqa}
\sigma(\mathcal{A})=
\begin{cases}
[0,\a]\cup[\b,\infty)&\text{ if }\a<\left({\pi\over 2L}\right)^2,\\
[0,\infty)&\text{otherwise},
\end{cases}
\end{gather}
where the number $\b$ is defined as follows. We denote by $\b(\mu)$ (here $\mu\in\mathbb{R}$) the smallest eigenvalue of the problem
\begin{gather*}
-u''=\lambda u\text{ in }\left(-L,0\right),\quad
u(-L)=0,\ u'(0)=\mu u(0).
\end{gather*}
It is straightforward to show that 
the function  $\mu\mapsto \b(\mu)$ is continuous, monotonically decreasing and moreover $\b(\mu)\underset{\mu\to -\infty}\to \left({\pi\over 2L}\right)^2$ and $\b(\mu)\underset{\mu\to +\infty}\to -\infty$. Whence, in particular, one can conclude that there exists one and only one point $\b$ satisfying
\begin{gather*}
\exists \mu<-\a r:\ \b=\b(\mu)={\a\mu \over \a r+\mu},
\end{gather*}
provided $\a<\left({\pi\over 2L}\right)^2$.\medskip

Now, we are in position to formulate the main results. 

\begin{theorem}\label{th1}
One has:
\begin{itemize}
 \item[(i)] Let the family $\left\{\lambda^{\eps}\in\sigma(\mathcal{A}\e)\right\}_{\eps\in\mathcal{E}}$ 
 have a convergent subsequence, i.e. $\lambda^{\eps}\to\lambda$ as $\eps=\eps'\to 0$. Then $\lambda\in \sigma(\mathcal{A})$.

 \item[(ii)] Let $\lambda\in \sigma(\mathcal{A})$. Then there exists a family $\left\{\lambda^{\eps}\in\sigma(\mathcal{A}\e)\right\}_{\eps\in\mathcal{E}}$ such that   $\liml_{\eps\to 0}\lambda\e=~\lambda$.

\end{itemize}

\end{theorem}

From \eqref{aqa} and Theorem \ref{th1} we immediately obtain the following 

\begin{corollary}
Let $\a<\left({\pi\over 2L}\right)^2$. Let $\delta>0$ be an arbitrary number satisfying $2\delta<\b-\a$. Then there exists $\eps_\delta>0$ such that
$$\sigma(\mathcal{A}\e)\cap (\a+\delta,\b-\delta)=\varnothing,\quad \sigma(\mathcal{A}\e)\cap (\a-\delta,\b+\delta)\not=\varnothing\quad\text{provided $\eps<\eps_\delta$.}$$

\end{corollary}

\section{Proof of Theorem \ref{th1}}
We present only the sketch of the proof since the main ideas are similar to the case of bounded domains $\Omega$ presented in \cite{CarKhrab}.\smallskip

Let $\left\{\lambda\e\in \sigma(\mathcal{A}\e)\right\}_{\eps\in\mathcal{E}}$ and $\lambda^{\eps}\to\lambda$ as $\eps=\eps'\to 0$. One has to show that $\lambda\in \sigma(\mathcal{A})$. In what follows we will use the index $\eps$ keeping in mind $\eps'$.

We denote
$${\Om}=(0,1)\times (-L,0),\quad 
{\Om}\e=\Omega\e\cap \left((0,1)\times \mathbb{R}\right),\quad \Ga=(0,1)\times\{0\}.$$
Recall  that $\eps^{-1}\in\mathbb{N}$, whence $\Omega\e+e_1=\Omega\e,$ where $e_1=(1,0),$
and thus $\mathcal{A}\e$ is a periodic operator with respect to the period cell $\Om\e$.

Using Floquet-Bloch theory (see, e.g, \cite{Kuchment}) one can represent the spectrum of $\mathcal{A}\e$ as a union of spectra of certain operators on $\Om\e$. We denote by $\widetilde{H}\e$ the space of functions from $L^2(\Om\e)$ and the scalar product defined by \eqref{He} with $\Om\e$ instead of $\Omega\e$.
Let us fix $\phi\in [0,2\pi)$. In  $\widetilde{H}\e$ we consider the sesquilinear form $\widetilde{\mathfrak{a}}^{\phi,\eps}$ defined by \eqref{ae} with $\Om\e$ instead of $\Omega\e$ and the definitional domain
$$\mathrm{dom}(\widetilde{\mathfrak{a}}^{\phi,\eps})=\left\{u\in H^1(\Om\e):\quad u(1,\cdot){=}e^{i\phi}u(0,\cdot)\right\}.$$
By $\widetilde{\mathcal{A}}^{\phi,\eps}$ we denote the operator associated with this form.
The spectrum of $\widetilde{\mathcal{A}}^{\phi,\eps}$ is purely discrete. We denote by $\{\widetilde\lambda^{\phi,\eps}_k\}_{k=1}^\infty$ the sequence of eigenvalues of $\widetilde{\mathcal{A}}^{\phi,\eps}$ arranged in ascending order and with account of their multiplicity. 
Then one has
\begin{gather}\label{floque}
\sigma(\widetilde{\mathcal{A}}\e)=\cupl_{k=1}^\infty I_k\e,\ \text{ where }I_k\e=\cupl_{\phi\in[0,2\pi)}\left\{\widetilde\lambda_{k}^{\phi,\eps}\right\}\text{ are compact intervals}.
\end{gather}

We also introduce the operator $\widetilde{\mathcal{A}}^{\phi}$ as the operator acting in $$\widetilde{H}=
\left\{U\in L^2(\Om)\oplus L^2(\Ga),\text{ the scalar product is defined by \eqref{Hlim}}
\text{ with }\Om,\Ga\text{ instead of }\Omega,\Gamma\right\}$$ and generated by the sesquilinear form $\widetilde{\mathfrak{a}}^\phi$ which is defined by \eqref{alim} (with $\Om,\Ga$ instead of $\Omega,\Gamma$) and definitional domain $\mathrm{dom}(\widetilde{\mathfrak{a}}^\phi)=\mathrm{dom}(\widetilde{\mathfrak{a}}^{\phi,\eps})\oplus L^2(\Ga)$.

\begin{lemma}\label{lmW2}
The spectrum of $\widetilde{\mathcal{A}}^{\phi}$ has the form
\begin{gather*}
\label{spectrum} \sigma(\widetilde{\mathcal{A}}^{\phi})=\{\a\}\cup
\{\widetilde\lambda_k^{\phi,-},k=1,2,3...\}\cup\{\widetilde\lambda_k^{\phi,+},k=1,2,3...\}.
\end{gather*}
The points $\widetilde\lambda_k^{\phi,\pm}, k=1,2,3...$ belong to the discrete
spectrum, $\a$ is a point of the essential spectrum and they are distributed  as follows:
$$0\leq\widetilde\lambda_1^{\phi,-}\leq \widetilde\lambda_2^{\phi,-}\leq ...\leq\widetilde\lambda_k^{\phi,-}\leq\dots\underset{k\to\infty}\to
\a< \widetilde\lambda_1^{\phi,+}\leq
\widetilde\lambda_2^{\phi,+}\leq ...\leq\widetilde\lambda_k^{\phi,+}\leq\dots\underset{k\to\infty}\to \infty.$$
Moreover if $\a<\left({\pi\over 2L}\right)^2$ then $\b< \widetilde\lambda_1^{\phi,+}.$
\end{lemma}
This lemma was proved in \cite{CarKhrab} for the case of Neumann boundary conditions on the lateral parts of $\partial\Om\e$. For the case of $\varphi$-periodic conditions the proof is similar. 
\smallskip

Now, in view of \eqref{floque} there exists $\phi\e\in [0,2\pi)$ such that $\lambda\e\in\sigma(\widetilde{\mathcal{A}}^{\phi\e\hspace{-1mm},\eps})$. We extract a convergent subsequence (for convenience still indexed by $\eps$):
\begin{gather}\label{phi}
\phi\e\to\phi\in [0,2\pi]\text{ as }\eps\to 0.
\end{gather}

Let $u\e$ be an eigenfunction of $\widetilde{\mathcal{A}}^{\phi\e\hspace{-1mm},\eps}$ corresponding to $\lambda\e$ with $\|u\e\|_{\widetilde{\mathcal{H}}\e}=1$.

We introduce the operator $\Pi\e:L^2(\cupl_{i=1}^{N(\eps)} B_i\e)\to L^2(\Gamma)$ defined as follows:
\begin{gather*}
\Pi\e u(x)= \suml_{i=1}^{N(\eps)}\left(|B_i\e|^{-1}\int_{B_i\e} u(\x)\d\x\right) \chi_{i}\e(x),
\end{gather*}
where $\chi_i\e$ is the characteristic function of the interval $\left[i\eps-\eps,i\eps\right]$.
Using the Cauchy inequality and  \eqref{qere+} one can easily obtain the estimate
\begin{gather}
\label{Pi_ineq} \|\Pi\e u\|^2_{L^2(\Gamma)}\leq \suml_{i=1}^{N(\eps)}{\varrho\e\eps|B_i\e|^{-1}}\int_{B_i\e}|u(\x)|^2(\varrho\e)^{-1}\d \x\leq C\|u\|^2_{\widetilde{\mathcal{H}}\e}.
\end{gather}
From \eqref{Pi_ineq} {and} $\|\nabla\e u\e\|^2_{L^2(\Omega\e)}=\lambda\e\leq C$,
we conclude that $\{u\e\}_{\eps\in\mathcal{E}}$ and $\{\Pi\e u\e\}_{\eps\in\mathcal{E}}$ are bounded in $H^1(\widetilde\Omega)$ and $L^2(\widetilde\Gamma)$, correspondingly. Then there is a subsequence (still indexed by $\eps$) and $u_1\in H^1(\Om)$, $u_2\in L^2(\Ga)$ such that
\begin{gather*}
u\e\rightharpoonup u_1\text{ in }H^1(\Om),\quad
\Pi\e u\e\rightharpoonup u_2\text{ in }L^2(\Ga).
\end{gather*} 
Also in view of the trace theorem and \eqref{phi} $u\e|_{\partial \widetilde\Omega}\to u_1$ in $L^2(\partial\widetilde\Omega)$, whence
$u_1(1,\cdot){=}e^{i\phi}u_1(0,\cdot)$, i.e. $U=(u_1,u_2)\in \mathrm{dom}(\widetilde{\mathfrak{a}}^\phi)$.

If $u_1=0$ then $\lambda=\a$, the proof is completely similar to the proof of this fact in \cite[Theorem 2.1]{CarKhrab}. Then in view of \eqref{aqa} $\lambda\in\sigma(\mathcal{A})$.\smallskip

Now, let $u_1\not= 0$. For an arbitrary $w\in \mathrm{dom}(\widetilde{\mathfrak{a}}^{\phi\e\hspace{-1mm},\eps})$ we have
\begin{gather}\label{in_eq_phi}
\int_{\Om\e}\nabla u\e(\x)\cdot \overline{\nabla w(\x)} \d \x=\lambda\e\int_{\Om\e} (\rho\e(\x))^{-1} u\e(\x) \overline{w(\x)} \d \x.
\end{gather}

Let $w_1\in C^{\infty}(\overline{\Om}),\ w_2\in C^{\infty}(\overline{\Ga})$, moreover $w_1(1,\cdot){=}e^{i\phi}w_1(0,\cdot).$
We set
$$w_1\e(\x)=w_1(\x)\left((e^{i(\phi\e-\phi)}-1)x_1+1\right),\ \x=(x_1,x_2).$$
It is easy to see that $w_1\e(\x)$ satisfies $w_1(1,\cdot){=}e^{i\phi\e}w_1(0,\cdot)$ and  \begin{gather}\label{appr1}
w_1\e\to w_1\text{ in }C^1(\overline{\Om})\text{ as }\eps\to 0.
\end{gather}

Using these functions we construct the test-function $w(x)$ by the formula 
\begin{gather*}
w(\x)= 
\begin{cases}
w_1\e(x)+\ds\suml_{i\in \I}(w_1\e(\x^{i,\eps})-w_1\e(\x))\varphi\left(\eps^{-1}|\x-\x^{i,\eps}|\right),& \x\in\Om,\\
\ds (h\e)^{-1}\left( w_2(\x^{i,\eps})-w_1\e(\x^{i,\eps})\right)x_2+w_1\e(\x^{i,\eps}),& \x=(x_1,x_2)\in T_i\e,\\
\ds{ w_2(\x^{i,\eps})},&\x\in B_i\e.
\end{cases}
\end{gather*}
Here $\x^{i,\eps}:=(i\eps,0)$, $\varphi\in C^\infty(\mathbb{R})$ satisfies $\varphi(t)=1$ as $t\leq {R\over 2}$ and $\varphi(t)=0$ as $t\geq {1\over 2}$, the constant $R\in (0,1)$ comes from \eqref{ass12}-\eqref{ass14}. It is clear that $w\e\in \mathrm{dom}(\mathfrak{a}^{\phi\e\hspace{-1mm},\eps})$.

We plug $w(\x)$ into \eqref{in_eq_phi} and pass to $\eps\to 0$.
Using the same arguments as in the proof of Theorem 2.1 from \cite{CarKhrab} (but with account of
 \eqref{appr1}) we obtain:
\begin{gather*} 
\int_\Omega {\nabla u_1}\cdot \overline{\nabla w_1}\d 
\x+\int_\Gamma \a r\left(u_1|_\Gamma-u_2\right)\overline{\left(w_1|_\Gamma-w_2\right)}\d x=\lambda \int_{\Om}u_1 \overline{w_1} \d \x+\lambda r\int_{\Ga}u_2 \overline{w_2} \d x.
\end{gather*}
By the density arguments this equality holds for an arbitrary $(w_1,w_2)\in\mathrm{dom}(\widetilde{\mathfrak{a}}^{\varphi})$ which implies
$\widetilde{\mathcal{A}} ^{\phi}U=\lambda U,\quad U=(u_1,u_2).$
Since $u_1\not=0$ then $\lambda\in\sigma(\widetilde{\mathcal{A}}^{\phi})$. 
But in view of \eqref{aqa} and Lemma \ref{lmW2} for each $\varphi$ one has  $\sigma(\widetilde{\mathcal{A}} ^{\phi})\subset \sigma(\mathcal{A} )$, therefore $\lambda\in\sigma(\mathcal{A})$.                                                                                                                                                                                                                                                                                                                                                                                                                                                      
The property (i) is proved.

\smallskip

The proof of the property (ii) repeats word-by-word the proof for bounded domains $\Omega$ presented in \cite{CarKhrab}. 

\section{References}

\renewcommand{\refname}{}    %%%% for this example
\vspace*{-26pt}              %%%% file only!

\end{document}